# Convexification of Queueing Formulas by Mixed-Integer Second-Order Cone Programming: An Application to a Discrete Location Problem with Congestion

**Amir Ahmadi-Javid**[1]

*Department of Industrial Engineering, Amirkabir University of Technology, Tehran, Iran*

ahmadi_javid@aut.ac.ir

**Pooya Hoseinpour**

*Desautels Faculty of Management, McGill University, Montreal, Quebec, Canada*

pooya.hoseinpour@mail.mcgill.ca

**Abstract.** Mixed-Integer Second-Order Cone Programs (MISOCPs) form a nice class of mixed-inter convex programs, which can be solved very efficiently due to the recent advances in optimization solvers. Our paper bridges the gap between modeling a class of optimization problems and using MISOCP solvers. It is shown how various performance metrics of M/G/1 queues can be molded by different MISOCPs. To motivate our method practically, it is first applied to a challenging stochastic location problem with congestion, which is broadly used to design socially optimal service networks. Four different MISOCPs are developed and compared on sets of benchmark test problems. The new formulations efficiently solve large-size test problems, which cannot be solved by the best existing method. Then, the general applicability of our method is shown for similar optimization problems that use queue-theoretic performance measures to address customer satisfaction and service quality.

---

[1] Corresponding author

# 1 Introduction

Second-order cone programming, also known as conic quadratic optimization, is a very promising branch of convex optimization from both theoretical and practical perspectives. Over recent years, different polynomial-time algorithms have been developed to solve Second-Order Cone Programs (SOCPs). The most popular solution methods for SOCPs are interior-point methods; which were initially proposed by Nemirovski and Scheinberg (1996), and then have been extended by many researchers because of their polynomial-time theoretical convergence and efficient computational performance; see Alizadeh and Goldfarb (2003) for a review. Ben-Tal and Nemirovski (2001a) also show that there is a polynomial-size polyhedral approximation for any second order cone, which pave the way to use linear programming to approximately solve SOCPs. For more technical details on these solution algorithms and new progresses, see, for example, Ben-Tal and Nemirovski (2001b), Boyd and Vandenberghe (2004), Benson and Saglam (2013), Kitahara and Tsuchiya (2016), Pena and Soheili (2017) and references therein. Most of these algorithms have also been implemented in powerful optimization solvers, which can solve large-scale problem instances in reasonable times. MATLAB toolboxes for optimization, such as SeDuMi and SDPT3, MOSEK, or CPLEX are some examples of well-known solvers for SOCPs.

Many classic decision problems can be modeled as SOCPs (Lobo, et al., 1998, Alizadeh and Goldfarb, 2003, and Nemirovski, 2006). Recent studies show that SOCPs can be used in various emerging application areas. Shivaswamy et al. (2006) present SOCPs for handling missing and uncertain data in machine learning. Maggioni et al. (2009) propose a two-stage stochastic SOCP in mobile ad-hoc networks. See and Sim (2010) consider a multi-period inventory control problem using robust optimization approach and find the policy by solving a tractable SOCP. Mak et al. (2014) consider an appointment scheduling problem and model it as an SOCP in a special case. Kocuk et al. (2016) consider an AC optimal power flow problem and propose three different relaxation models in format of SOCPs. Coutinho et al. (2016) apply SOCPs in a branch-and-bound algorithm to solve close-enough traveling salesman problem. Vu et al. (2016) developed SOCPs for wireless communications design. Chen and Zhu (2017) propose an SOCP model for a two-stage network data envelopment analysis.

The promising achievements in solving SOCPs and Linear Programs (LPs) have recently been used to make Mixed-Integer SOCP (MISOCP) as a tractable class of mixed-integer convex programming to solve nice discrete decision problems in industrial sizes. There are two major groups of algorithms for solving MISOCPs. The first one is based on branch-and-bound (B&B) method where solving an SOCP in each node is required. The second one is to use outer-approximation branch-and-cut (B&C) method. CPLEX can apply both approaches, as well as an additional option where the best approach is cleverly selected. For more details and comparisons see Drewes, 2009, Drewes and Ulbrich (2012), and Bonami and Tramontani (2015). Moreover, several kinds of cuts have been developed for MISOCPs; many of them are applied in the existing solvers to empower them to efficiently solve large-scale MISOCPs (Atamtürk & Narayanan, 2010, 2011; Bonami, 2011; Dadush, 2011; Andersen & Jensen, 2013; Belotti et al., 2013, 2015; Goez, 2013; Kılınç-Karzan & Yıldız, 2015; Kılınç-Karzan, 2015; Modaresi et al., 2015, 2016; Pattanayak & Narayanan, 2017; Kılınç et al., 2017).

Based on the recent progress on developing solution algorithms for MISOCPs, nowadays these optimization problems are efficiently solvable in many practical cases. Recently, MISOCPs have been used to model and solve many challenging applied problems. Some of the interesting applications, reviewed by Benson and Saglam (2013), are portfolio optimization (Vielma et al., 2008; Bonami & Lejeune, 2009), option pricing (Pinar, 2013), telecommunication network design (Hijazi et al., 2013), transmission in cellular networks (Cheng et al., 2012), power distribution (Taylor & Hover, 2012), battery swapping stations on freeway networks (Mak et al., 2013), Euclidean k-centering (Brandenberg & Roth, 2009), location-inventory planning (Ahmadi-Javid & Azad, 2010; Atamtürk et al., 2012), and scheduling and logistics (Du et al., 2011). More recently, Pinar (2013) develop an MISOCP for lower hedging of American contingent claims. Miyashiro and Takano (2015) propose an MISOCP for explanatory variable selection in a multiple linear regression model. Borraz-Sánchez et al. (2016) present an MISOCP relaxation for a gas expansion-planning problem under steady-state conditions. Han et al. (2016) use an MISOCP to develop an approximation algorithm for optimal learning in linear regression. He et al. (2017) propose an MISOCP as an approximate model for designing urban electronic-car sharing system. Tanaka and Kobayashi (2017) present an MISOCP for optimal fuel route problem in which the shipping route and its speed are found to minimize the total fuel consumption between two ports.

This paper addresses a class of discrete optimization problems that deal with performance metrics of M/G/1 queues. These metrics are the expected waiting times in the queue and system, or the total congestion in the queue and system. They may appear in the objective functions or constraints to control the congestion level, or provide a specific service level (Boffey et al. 2007). We show that important metrics of an M/G/1 queue can be represented by different SOCP formulations in a flexible setting.

To simply present our idea and show its practical advantages, we start by examining a traditional and challenging location problem, studied by several papers for three decades. We present different MISOCPs for this location problem, and then we compare the MISOCPs with each other, as well as the best known exact solution method for this problem. Then, we discuss how our reformulation idea can be applied in a general setting.

The remainder of the paper is organized as follows. Section 2 provides some required backgrounds on queue systems and SOCPs. Section 3 introduces the classic location problem with congestion and provides the basic non-convex integer programing model of this problem. Section 4 presents different MISOCP reformulations for this problem, and Section 5 carries out a comprehensive numerical study to assess the performance of these MISOCPs. Section 6 demonstrates how our method can be used in other similar problems. Section 7 provides an outlook.

## 2  Preliminaries

This section presents some preliminaries needed throughout the paper. Section 2.1 presents an overview on M/G/1 queue systems. Section 0 provides a brief introduction to MISOCPs. Section 2.3 discusses those forms of MISOCPs that are acceptable by CPLEX.

### 2.1 An overview on M/G/1 queue systems

Using Kendall's notation (Kendall, 1953), an M/G/1 queue system is a single server queue system, which can be modeled as a continuous-time stochastic process. The arrival of customers occurs based on a Poisson process with intensity $\lambda$. The service times $S_1, S_2, \ldots \sim S$ of the customers are independent and generally distributed with finite mean $E(S)$ and variance $\text{var}(S)$; the departure rate is denoted by $\mu = 1/E(S)$. For stability of the queue status, it is assumed that $\rho := \lambda/\mu = \lambda E(S) < 1$.

The stationary queue-length distribution is equal to the distribution of the number of customers in the system on departure instants, for which, using Pollaczek–Khinchine formula (Gross, 2008), the probability generating function can be calculated as follows:

$$P(z) = \frac{(1-\rho)LS_S(\lambda - \lambda z)(1-z)}{LS_S(\lambda - \lambda z) - z} \qquad |z| \leq 1, \tag{1}$$

where $LS_S$ is the Laplace-Stieltjes transform of the distribution of the service time $S$. Using (1), the expected numbers of people in the queue and in the system; denoted by $L^q$ and $L$, respectively, can be calculated as

$$L^q = \frac{\rho^2 + \lambda^2 \sigma^2}{2(1-\rho)} \tag{2}$$

$$L = \rho + L^q = \rho + \frac{\rho^2 + \lambda^2 \sigma^2}{2(1-\rho)}, \tag{3}$$

where $\sigma^2 = \text{var}(S)$ is the variance of the service-time distribution. Hence, the mean times that a customer spends in the queue and in the system; denoted by $W^q$ and $W$, respectively, are given by

$$W^q = \frac{\rho + \lambda \mu \sigma^2}{2(\mu - \lambda)} \tag{4}$$

$$W = W^q + \mu^{-1} = \frac{\rho + \lambda \mu \sigma^2}{2(\mu - \lambda)} + \mu^{-1}, \tag{5}$$

which follow from (2), (3), and Little's law, stating that $L^q = \lambda W^q$ and $L = \lambda W$. For more details on computing the waiting-time distributions of M/G/1 queues and other recent related developments, see e.g., Shortle et al. (2004, 2007), Connor & Kendall (2015), Sigman, K. (2016), Baron & Kerner (2016), Li et al. (2017), and references therein.

One should note that the expected total waiting times in the queue and in the system, denoted by $TW^q$ and $TW$, respectively, are nothing but $L^q$ and $L$, respectively. These metrics are used to express the overall congestion in a service center. Hence, we have

$$TW^q = \lambda W^q = L^q \tag{6}$$

$$TW = \lambda W = L. \tag{7}$$

These metrics are usually used to express the overall congestion in a service center.

## 2.2 Mixed-integer second-order cone programming

Let $\kappa \subseteq \mathbb{R}^m$ be a proper cone, i.e., a pointed and closed convex cone with a nonempty interior. Then, this cone induces a partial ordering on $\mathbb{R}^m$, denoted by $\succcurlyeq_\kappa$, where we have

$$a \succcurlyeq_\kappa b \Leftrightarrow a - b \succcurlyeq_\kappa 0 \Leftrightarrow a - b \in \kappa.$$

The optimization problem

$$\min_{x \in \mathbb{R}^n} c^T x$$
$$A_i x - b_i \succcurlyeq_{\kappa_i} 0, i = 1, \dots, p$$
$$x \in P,$$

is called a *conic program* where $c \in \mathbb{R}^n, b \in \mathbb{R}^n, A_i$ is an $m_i \times n$ matrix, $\kappa_i \subseteq \mathbb{R}^{m_i}, i = 1, \dots, p$ are proper cones, and $P$ is a polyhedron in $\mathbb{R}^n$. Actually, a conic program is a linear optimization problem with generalized linear inequalities

$$A_i x - b_i \succcurlyeq_{\kappa_i} 0, i = 1, \dots, p.$$

A *second-order cone* is a proper cone defined by

$$S^m = \left\{ z \in \mathbb{R}^m : z_m \geq \sqrt{\sum_{i=1}^{m-1} z_i^2} \right\},$$

for $m \geq 2$. This cone is also called a *Lorentz, quadratic,* or *ice-cream* cone. An SOCP (Second-Order Conic program), also called a Conic Quadratic Program (CQP), is a conic program for which all the cones $\kappa_i, i = 1, \dots, p$ are Second-Order Cones (SOCs). Therefore, by defining $[A_i, b_i] = \begin{bmatrix} E_i & e_i \\ \beta_i^T & \delta_i \end{bmatrix}$ the above **conic form** of an SOCP can simply be rewritten as follows:

$$\min_{x \in \mathbb{R}^n} c^T x$$
$$\|E_i x - e_i\| \leq \beta_i^T x - \delta_i \quad i = 1, \cdots, p$$
$$x \in P.$$

This form is here called the **primary form** of an SOCP. When all $\beta_i$ and $\delta_i, i = 1, \cdots, p$, are zeroes, the SOCP reduces to an LP. If all $\beta_i, i = 1, \cdots, p,$ are zeroes the SCOP becomes a convex quadratically-constrained linear program. SOCPs are important convex optimization problems which are polynomially solvable.

Each constraint $\|E_i x - e_i\| \leq \beta_i^T x - \delta_i$ is called an *SOC constraint*. A set of constraints (or any subset of $\mathbb{R}^n$) is called *SOC representable* if it can equivalently be expressed by a finite set of SOC constraints. A mathematical program is called an MISOCP (Mixed-Integer SOCPs) when some of decision variables $x_i$, $i = 1, \cdots, n$ are restricted to be integer.

## 2.3 MISOCP forms acceptable by CPLEX

Using CPLEX, we are able to exactly solve mixed-integer programs with quadratic constraints, in the form of $1/2\, x^T Q x + \beta^T x + \alpha \leq 0$, where $Q$ must be checked to be a positive semi-definite matrix, i.e., $Q \succcurlyeq 0$. Moreover, to handle MISOCPs, CPLEX accepts the following two forms of constraints for a given semi-definite matrix $Q \succcurlyeq 0$:

I)     $x^T Q x \leq y^2$ where $y \geq 0$

II)    $x^T Q x \leq yz$ where $y, z \geq 0$.

The variable $y$ (or $z$) can be replaced by a positive affine transformation of some non-negative variables, i.e., $y = \sum_{i \in G} c_i w_i$ where $w_i \geq 0$, $i \in G$, are decision variables, and $c_i > 0$, $i \in G$, are positive constants. However, carefully note that the variable $y$ (or $z$) cannot be replaced by an affine transformation that is not independently non-negative, even if it always non-negative by considering the other constraints of the problem. Also, note that

$$\sum_{i=1}^{n} x_i^2 \leq y^2$$

with $y \geq 0$ is a special case of the first from with $Q = I$.

Our computational experiment indicate that CPLEX can handle the form-I constraints more efficiently. Thus, a form-II constraint is suggested to be replaced by the following alternative

$$2x^T Q x + y^2 + z^2 \leq (y + z)^2 \text{ where } y, z \geq 0,$$

which is a form-I constraint because $y + z$ is always non-negative for any $y, z \geq 0$. See IBM (2017) for more details.

When all SOC constraints of an SOCP are transformed to form-I and form-II constraints, the resulting formulation is here called the **secondary form** of the SOCP.

The above points are very important to solve MISOCPs by a solver such as CPLEX. In fact, one should make some further changes to a given MISOCP to make it acceptable for an MISOCP solver. This consideration results in introducing new variables and constraints. In Section 4, it will be seen that MISOCPs with the same primary-form formulations but different secondary-form formulations may have very dissimilar performance.

## 3 Stochastic location problem with congestion

This section introduces a congested location-allocation problem considered by many studies under differences that slightly affect the structures of the resulting optimization problems (refer to the comprehensive review by Berman and Krass, 2015). The problem establishes facilities, determines the capacities of the established facilities, and determines the established facility that provides service to each customer (demand zone) in order to minimize a performance metric that compromises between the cost of the service network and the clients' congestion and traveling costs. Each customer makes a demand stream evolving according to a Poisson process, and each facility can be molded as an M/G/1 queue system (see Section 2.1).

To have a formal definition of the problem, first let us define our notation. Suppose that $I, J$, and $K$ denote the sets of available facilities, customers, and service-capacity levels, respectively. The following parameters and decision variables are considered throughout the modeling:

**Parameters:**

$f_{ik}$   The fixed cost for establishing facility $i \in I$ with service-capacity level $k \in K$.

$d_{ij}$   The traveling cost from the location of customer $j \in J$ to facility $i \in I$ per service.

$\lambda_j$   The demand rate of customer $j \in J$.

$\mu_{ik}$   The service rate facility $i \in I$ with service-capacity level $k \in K$.

$\sigma_{ik}$   The service-time standard deviation at facility $i \in I$ with service-capacity level $k \in K$.

$w_i$   The congestion (waiting) cost per unit time at facility $i \in I$.

**Decision variables:**

$x_{ik}$   A binary variable that takes 1 if facility $i \in I$ is established with service-capacity level

$k \in K$, and 0 otherwise.

$y_{ikj}$  A binary variable that takes 1 if customer $j \in J$ is served by facility $i \in I$ with service-capacity level $k \in K$, and 0 otherwise.

Next, our service network design problem can be formulated as follows:

$$\min \sum_{i \in I} \sum_{k \in K} f_{ik} x_{ik} + \sum_{i \in I} \sum_{k \in K} w_i TW_{ik}\left(\sum_{j \in J} \lambda_j y_{ikj}, \mu_{ik}, \sigma_{ik}\right) + \sum_{i \in I} \sum_{k \in K} \sum_{j \in J} d_{ij} \lambda_j y_{ikj} \quad (8)$$

s.t.

$$\sum_{i \in I} \sum_{k \in K} y_{ikj} = 1 \qquad j \in J \quad (9)$$

$$y_{ikj} \leq x_{ik} \qquad i \in I, k \in K, j \in J \quad (10)$$

$$\sum_{k \in K} x_{ik} \leq 1 \qquad i \in I \quad (11)$$

$$\sum_{j \in J} \lambda_j y_{ikj} \leq \mu_{ik} \qquad i \in I, k \in K \quad (12)$$

$$x_{ik} \in \{0,1\} \qquad i \in I, k \in K \quad (13)$$

$$y_{ikj} \in \{0,1\} \qquad i \in I, k \in K, j \in J \quad (14)$$

where $TW_{ik}(\cdot)$ is the total waiting at facility $i \in I$ with service-capacity level $k \in K$, which by (3) becomes as follows:

$$TW_{ik}\left(\sum_{j \in J} \lambda_j y_{ikj}, \mu_{ik}, \sigma_{ik}\right) = \frac{\left(\sum_{j \in J} \lambda_j y_{ikj}\right)^2 (1 + \mu_{ik}^2 \sigma_{ik}^2)}{2\mu_{ik}\left(\mu_{ik} - \sum_{j \in J} \lambda_j y_{ikj}\right)} + \frac{\sum_{j \in J} \lambda_j y_{ikj}}{\mu_{ik}}. \quad (15)$$

The first term in objective function (8) is the annual cost of establishing facilities. The second term is the yearly total congestion cost at facilities, while the third term is the yearly total accessing cost of customers. Constraint set (9) ensures that each customer is assigned to only one established facility. Constraint set (10) forbids assigning a customer to a non-established facility. Constraint set (11) forces that at most one service-capacity level is selected for each established facility. Constraint set (12) ensures the required steady-state condition at each established facility. Constraint sets (13) and (14) guarantee that all decision variables are binary.

The model (8)-(15) is an integer non-convex minimization problem, which cannot be optimally solved in practical scales using existing general-purpose solvers. The problem is recently solved in Vidyarthi and Jayaswal (2014) by an efficient exact constraint-generation method, which is an extension of the method used by Elhedhli (2006) for the M/M/1 case. In fact, before 2006, only Lagrangian heuristics were available to tackle this problem for the special case where facilities are represented by M/M/1 queues (Berman and Krass, 2015).

In the next section, it is shown that the model (8)-(15) can be cast as different MISOCPs, which are efficiently solvable using MISOPC-solvers such as CPLEX. In Section 5, their performance is assessed and compared with the best exiting algorithm given by Vidyarthi and Jayaswal (2014).

## 4 MISOCP reformulations for location problem

The first part of this section provides four secondary-form MISOCP reformulations for the model (8)-(15), as well as their corresponding primary-form formulations. The second part structurally compares the formulations.

The nonlinear terms in (4) that represent the congestion costs are the only complicating factors. Thus, in the sequel, our focus is to demonstrate how these nonlinear terms are SOC representable and how they can be rewritten by secondary-form MISOCP reformulations, which can be solved using existing solvers such as CPLEX (see Section 2.3). In each subsection, the associated primary-form formulation is also presented.

### 4.1 The first MISOCP reformulation

The model (8)-(14) can be rewritten as the following secondary-form MISOPC:

$$\min \sum_{i \in I} \sum_{k \in K} f_{ik} x_{ik} + \sum_{i \in I} \sum_{k \in K} \frac{w_i (1 + \mu_{ik}^2 \sigma_{ik}^2) r_{ik}}{2 \mu_{ik}} + \frac{w_i \sum_{j \in J} \lambda_j y_{ikj}}{\mu_{ik}} + \sum_{i \in I} \sum_{k \in K} \sum_{j \in J} d_{ij} \lambda_j y_{ikj} \quad (16)$$

s.t.

$$2 \left( \sum_{j \in J} \lambda_j y_{ikj} \right)^2 + r_{ik}^2 + t_{ik}^2 \leq (r_{ik} + t_{ik})^2 \qquad i \in I, k \in K \quad (17)$$

$$t_{ik} = \mu_{ik} - \sum_{j \in J} \lambda_j y_{ikj} \qquad i \in I, k \in K \qquad (18)$$

$$t_{ik} \geq 0 \qquad i \in I, k \in K \qquad (19)$$

$$r_{ik} \geq 0 \qquad i \in I, k \in K \qquad (20)$$

(9)-(14).

To show this, for each $i \in I, k \in K$ let us introduce the auxiliary variables $r_{ik}$ and the new constraint $r_{ik} \geq \frac{(\sum_{j \in J} \lambda_j y_{ikj})^2}{\mu_{ik} - \sum_{j \in J} \lambda_j y_{ikj}}$. Using these, the model (8)-(14) can be transformed to

$$\min \sum_{i \in I} \sum_{k \in K} f_{ik} x_{ik} + \sum_{i \in I} \sum_{k \in K} \frac{w_i (1 + \mu_{ik}^2 \sigma_{ik}^2) r_{ik}}{2 \mu_{ik}} + \frac{w_i \sum_{j \in J} \lambda_j y_{ikj}}{\mu_{ik}} + \sum_{i \in I} \sum_{k \in K} \sum_{j \in J} d_{ij} \lambda_j y_{ikj} \qquad (16)$$

s.t.

$$\left( \sum_{j \in J} \lambda_j y_{ikj} \right)^2 \leq r_{ik} \left( \mu_{ik} - \sum_{j \in J} \lambda_j y_{ikj} \right) \qquad i \in I, k \in K \qquad (21)$$

(9)-(14), (20).

Constraints (21) are hyperbolic constraints, and can be transformed to the following form-II constraints:

$$\left( \sum_{j \in J} \lambda_j y_{ikj} \right)^2 \leq r_{ik} t_{ik} \qquad i \in I, k \in K$$

if constraints $t_{ik} = \mu_{ik} - \sum_{j \in J} \lambda_j y_{ikj}$, $i \in I, k \in K$, are added, where variables $t_{ik}$ are non-negative by (12). As it is stated in Section 2.3, to improve computational performance, it is better to replace these constraints by form-I constraints. Therefore, to make (21) acceptable for MISOCP solvers, they can be rewritten as (17)-(19). The equivalent primary-form SOC constraints are given by

$$\left\| \begin{array}{c} 2 \sum_{j \in J} \lambda_j y_{ikj} \\ r_{ik} - \mu_{ik} + \sum_{j \in J} \lambda_j y_{ikj} \end{array} \right\| \leq r_{ik} + \mu_{ik} - \sum_{j \in J} \lambda_j y_{ikj} \qquad i \in I, k \in K.$$

## 4.2 The second MISOCP reformulation

The model (8)-(14) can be transformed to the following secondary-form MISOCP formulation:

$$\min \sum_{i \in I} \sum_{k \in K} f_{ik} x_{ik} + \sum_{i \in I} \sum_{k \in K} \frac{w_i(1 + \mu_{ik}^2 \sigma_{ik}^2) s_{ik}}{2} + \frac{w_i(1 - \mu_{ik}^2 \sigma_{ik}^2) \sum_{j \in J} \lambda_j y_{ikj}}{2\mu_{ik}}$$
$$+ \sum_{i \in I} \sum_{k \in K} \sum_{j \in J} t_{ij} \lambda_j y_{ikj} \qquad (22)$$

s.t.

$$2 \sum_{j \in J} \lambda_j y_{ikj}^2 + s_{ik}^2 + t_{ik}^2 \leq (s_{ik} + t_{it})^2 \qquad i \in I, k \in K \qquad (23)$$

$$s_{ik} \geq 0 \qquad (24)$$

(9)-(14), (18)-(19).

To accept this, new constraints $s_{ik} \geq \frac{\sum_{j \in J} \lambda_j y_{ikj}}{\mu_{ik} - \sum_{j \in J} \lambda_j y_{ikj}}$, $i \in I, k \in K$, should be considered to obtain

$$\min \sum_{i \in I} \sum_{k \in K} f_{ik} x_{ik} + \sum_{i \in I} \sum_{k \in K} \frac{w_i(1 + \mu_{ik}^2 \sigma_{ik}^2) s_{ik}}{2} + \frac{w_i(1 - \mu_{ik}^2 \sigma_{ik}^2) \sum_{j \in J} \lambda_j y_{ikj}}{2\mu_{ik}}$$
$$+ \sum_{i \in I} \sum_{k \in K} \sum_{j \in J} d_{ij} \lambda_j y_{ikj} \qquad (22)$$

s.t.

$$\sum_{j \in J} \lambda_j y_{ikj} \leq s_{ik} \left( \mu_{ik} - \sum_{j \in J} \lambda_j y_{ikj} \right) \qquad i \in I, k \in K \qquad (25)$$

(9)-(14), (24).

Then, Constraints (25) can be transformed to the constraints below considering the fact that variables $y_{ikj}$ $i \in I, k \in K, j \in J$, are 0-1 valued and can be replaced by $y_{ikj}^2$. Therefore, one can use the following set of constraints instead of (25):

$$\sum_{j \in J} \lambda_j y_{ikj}^2 \leq s_{ik} \left( \mu_{ik} - \sum_{j \in J} \lambda_j y_{ikj} \right) \qquad i \in I, k \in K,$$

which are equivalent to (18), (19), (23). The corresponding primary-form SOC constraints are given by

$$\left\| \begin{matrix} \sqrt{2\lambda_1} y_{i1k} \\ \vdots \\ \sqrt{2\lambda_{|J|}} y_{i|J|k} \\ r_{ik} - \mu_{ik} + \sum_{j \in J} \lambda_j y_{ikj} \end{matrix} \right\| \leq r_{ik} + \mu_{ik} - \sum_{j \in J} \lambda_j y_{ikj} \qquad i \in I, k \in K.$$

Note that this formulation is not applicable for the case that variables $y_{ikj}$ are not binary.

### 4.3 The third MISOCP reformulation

The third secondary-form MISOCP formulation for the model (8)-(14) can be given by

$$\min \sum_{i \in I} \sum_{k \in K} f_{ik} x_{ik} + \sum_{i \in I} \sum_{k \in K} \frac{w_i (1 + \mu_{ik}^2 \sigma_{ik}^2) s_{ik}}{2} + \frac{w_i (1 - \mu_{ik}^2 \sigma_{ik}^2) \sum_{j \in J} \lambda_j y_{ikj}}{2 \mu_{ik}} \\ + \sum_{i \in I} \sum_{k \in K} \sum_{j \in J} d_{ij} \lambda_j y_{ikj} \qquad (22)$$

s.t.

$$4 \left( \sum_{j \in J} \lambda_j y_{ikj} \right)^2 + p_{ik}^2 \leq q_{ik}^2 \qquad i \in I, k \in K. \qquad (26)$$

$$p_{ik} = s_{ik} \mu_{ik} - \sum_{j \in J} \lambda_j y_{ikj} - \mu_{ik} + \sum_{j \in J} \lambda_j y_{ikj} \qquad i \in I, k \in K \qquad (27)$$

$$q_{ik} = s_{ik} \mu_{ik} - \sum_{j \in J} \lambda_j y_{ikj} + \mu_{ik} - \sum_{j \in J} \lambda_j y_{ikj} \qquad i \in I, k \in K \qquad (28)$$

$$p_{ik}, q_{ik} \geq 0 \qquad i \in I, k \in K \qquad (29)$$

(9)-(14), (24).

To obtain this formulation, by a direct calculation one can see that (25) is representable as

$$\left( \sum_{j \in J} \lambda_j y_{ikj} \right)^2 \leq \left( s_{ik} \mu_{ik} - \sum_{j \in J} \lambda_j y_{ikj} \right) \left( \mu_{ik} - \sum_{j \in J} \lambda_j y_{ikj} \right) \qquad i \in I, k \in K, \qquad (30)$$

which are hyperbolic constraints. Salimian and Gürel (2013) use the above representation to model a problem in a make-to-order supply chain with cross-docking terminals studied by Vidyarthi et al. (2009). They rewrote the above constraints as (26)-(29) to make them acceptable for CPLEX. Note that here we have improved their secondary-form formulation by eliminating unnecessary variables and linear constraints that they used to present the first term in (26). The constraints (26)-(29) can be represented by the following primary-form SOC constraints:

$$\left\| \begin{array}{c} 2\sum_{j \in J} \lambda_j y_{ikj} \\ s_{ik}\mu_{ik} - \sum_{j \in J} \lambda_j y_{ikj} - \mu_{ik} + \sum_{j \in J} \lambda_j y_{ikj} \end{array} \right\| \leq s_{ik}\mu_{ik} - \sum_{j \in J} \lambda_j y_{ikj} + \mu_{ik} - \sum_{j \in J} \lambda_j y_{ikj} \qquad i \in I, k \in K.$$

### 4.4 The fourth MISOCP reformulation

The fourth secondary-form MISOCP for the model (8)-(14) can be given by

$$\min \sum_{i \in I} \sum_{k \in K} f_{ik} x_{ik} + \sum_{i \in I} \sum_{k \in K} \frac{w_i(1 + \mu_{ik}^2 \sigma_{ik}^2) s_{ik}}{2} + \frac{w_i(1 - \mu_{ik}^2 \sigma_{ik}^2) \sum_{j \in J} \lambda_j y_{ikj}}{2\mu_{ik}}$$
$$+ \sum_{i \in I} \sum_{k \in K} \sum_{j \in J} d_{ij} \lambda_j y_{ikj} \qquad (22)$$

s.t.

$$2 \left( \sum_{j \in J} \lambda_j y_{ikj} \right)^2 + t_{ik}^2 + v_{ik}^2 \leq (t_{ik} + v_{it})^2 \qquad i \in I, k \in K. \qquad (31)$$

$$v_{ik} = s_{ik}\mu_{ik} - \sum_{j \in J} \lambda_j y_{ikj} \qquad i \in I, k \in K \qquad (32)$$

$$v_{ik} \geq 0 \qquad i \in I, k \in K \qquad (33)$$

(9)-(14), (18)-(19), (24).

This is directly obtained by reformulating form-II Constraints (30) as from-I Constraints (31) using the method explained in Section 2.3. The primary-form formulation of this MISOCP is the same as the one presented for the third secondary-form MISOCP given in Section 4.1.3.

### 4.5 Comparison of MISOCPs

Considering that the proposed MISOCP formulations are here solved using CPLEX, their secondary-form formulations should be compared. The characteristics of these formulations are summarized in Table 1. From this table, it may be expected that 1-MISOCP and 2-MISOCP have better computational performance than 3-MISOCP and 4-MISOCP because they have less real additional variables and less additional constraints. However, note that CPLEX solver is empowered by a pre-solving step, in which the models are first intelligently examined for reduction opportunities before solving; such as ignoring repeated or redundant constraints. In addition, probing techniques is applied; e.g., tentatively fixing some binary variable (Savelsbergh, 1994). Our computational results show that the integer-relaxation bounds of 1-MISOCP and 4-MISOCP are the same, and they are better (larger) than those of 2-MISOCP and 3-MISOCP. This can help 1-MISOCP and 4-MISOCP to reach small gaps much faster. Our computational results indicate that 1-MISOCP is much more efficient than 4-MISOCP in large scales. Moreover, the integer-relaxation bound of 4-MISOCP is often better than 3-MISOCP.

Table 1. Comparison of proposed MISOCPs.

| Secondary-form MISOCP | Objective | Constraints | Number of additional real variables | Number of additional constraints | Number of form-I constraints |
|---|---|---|---|---|---|
| 1-MISOCP | (16) | (9)-(14), (17)-(20) | $2\|I\| \times \|K\|$ | $\|I\| \times \|K\|$ | $\|I\| \times \|K\|$ |
| 2-MISOCP | (22) | (9)-(14), (18)-(19), (23)-(24) | $2\|I\| \times \|K\|$ | $\|I\| \times \|K\|$ | $\|I\| \times \|K\|$ |
| 3-MISOCP | (22) | (9)-(14), (24), (26)-(29) | $3\|I\| \times \|K\|$ | $2\|I\| \times \|K\|$ | $\|I\| \times \|K\|$ |
| 4-MISOCP | (22) | (9)-(14), (18)-(19), (24), (31)-(33) | $3\|I\| \times \|K\|$ | $2\|I\| \times \|K\|$ | $\|I\| \times \|K\|$ |

## 5  Computational results

In this section, the proposed MISOCPs are solved by IBM ILOG CPLEX Optimization Studio 12.6.1. We use a PC with a dual-core 2.9 GHz processor and 30GB RAM, operating Windows 7, 64-bit.

Table 2 compares the four secondary-form MISOCPs presented in Section **Error! Reference source not found.** on large-size test problems used in Vidyarthi and Jayaswal (2014), with 25 facilities, 5 service-capacity levels, and 400 customers. To solve MISOCPs, we allow the CPLEX solver to intelligently choose between SOCP B&B and outer-approximation B&C by setting parameter miqcpstrat = 0. This table shows that 1-MISOCP and 4-MISOCP perform better than 2-MISOCP and 3-MISOCP in finding optimal solutions in less run times as well as less memory usage. One can also see that the integer-relaxation bounds of 1-MISOCP and 4-MISOCP are always the same, and are better than those of the other two models. The integer-relaxation bound of 3-MISOCP is always better than that of 3-MISOCP expecting for test problem I12.

For the same instances given in Table 2, Table 3 compares the run times of solving our two selected models 1-MISOCP and 4-MISOCP with those required by the cutting-plane algorithm given by Vidyarthi and Jayaswal (2014) to achieve optimality-gap limit 0.001%. Table 3 also reports the total cost and the percentages of its components for each instance to make sure that the instances have sensibly balanced objective functions. As shown, our proposed models are successful in reaching the exact optimal solution in reasonable run times for all instances, while the cutting-plane method cannot optimally solve the instance I12 within 3 run-time hours. In four instances I02, I05, I06, and I07, the cutting-plane method has smaller run times compared to the two MISOCPs. To see what will happen in a larger scale, we test three alternatives on larger sizes.

Table 4 considers test problems that are initially used by Holmberg et al. (1999) for a capacitated facility location problem. Medium-size instances of these test problems were used later by Elhedhli (2006) for designing a service system where the service-capacity levels are set to be 3. Here, we consider the instances of these test problems with the largest sizes, and with 10 service-capacity levels. First, we generate all of our parameters in the same way considered by Elhedhli (2006). In order to generate the fixed cost of establishing facilities, $f_{ik}$, we similarly use the following formula to reflect the economy of scale:

$$f_{ik} = \left(\frac{f_i}{b_i}\right)^{\frac{|K|-1}{|K|-1+k-1}} \times \mu_{ik},$$

where $f_i, b_i$ are the facility cost and capacity in the original test problems given in Holmberg et al. (1999), respectively. Moreover, as our queue systems at facilities are M/G/1 instead of M/M/1, which is used in Elhedhli (2006), the new parameters $\sigma_{ik}$ are needed to be generated. They are randomly generated from the interval $[1/\mu_{ik}, 3/\mu_{ik}]$.

Table 4 summarizes the computational results on above-mentioned test problems generated based on Holmberg et al. (1999). In most of the cases, 1-MISOCP and 4-MISOCP reach to the optimality very efficiently, even within a few minutes for some test problems. However, the cutting-plane algorithm fails to reach its pre-defined optimality gap 0.001% in 3 hours for most of the test problems (expecting for three instances P64, P68, and P69). For test problems P57, P58, P62, and P66, none of the two MISOCPs and the cutting-plane method can find the optimal solutions within 3 hours, but the optimality gaps of the solutions found by the two MISOCPs are significantly better than those provided by the cutting-plane algorithm. This numerical analysis clearly indicates that our selected MISOCPs are both efficient and stable in solving large-size instances and significantly outperform the cutting-plane method. In other words, on the test problems of Table 4, the cutting plane method is completely dominated, while 1-MISOCP competes with 4-MISOCP.

To determine the best formulations between 1-MISOCP and 4-MISOCP, Table 5 carries out an additional analysis on test problems with larger sizes. These test problems are generated by combining each pair of consecutive test problems given in Table 4, which results in instances with 60 facilities, 10 service-capacity levels, and 400 customers. The comparison indicates that 1-MISOCP outperforms 4-MISOCP in all cases under each one of the two different run-time limits. Hence, 1-MISOCP seems to be the best candidate for solving the proposed congested location problem in a large scale.

# 6 MISCOPs for queueing formulas in a general setting

This section demonstrates our results in a general setting such that they can be used in other problems involving one of the performance metrics of an M/G/1 queue system. Let us consider the case that the service rate $\mu$, the variance of the service time $\sigma^2$, the arrival rate $\lambda$ can be represented by

$$\mu = \sum_{p \in P} \mu_p X_p$$

$$\sigma^2 = \sum_{p \in P} \sigma_p^2 X_p$$

$$\lambda = \sum_{p \in P} \sum_{r \in R} \lambda_r Y_{p,r}$$

where $\mu_p$, $\sigma_p^2$, and $\lambda_r$ are given constants and $X_p$, $p \in P$ and $Y_{p,r}$, $p \in P, r \in R$ are decision variables that satisfy the following conditions:

$$\sum_{p \in P} \sum_{r \in R} \lambda_r Y_{p,r} \leq \sum_{p \in P} \mu_p X_p \qquad (34)$$

$$\sum_{p \in P} X_p = 1 \qquad (35)$$

$$X_p \in \{0,1\} \qquad p \in P \qquad (36)$$

$$0 \leq Y_{p,r} \leq U_{p,r} X_p \qquad p \in P, r \in R \qquad (37)$$

for some given constant upper bounds $U_{p,r}$, $p \in P, r \in R$. Then, we show how the constraint

$$PM(X,Y) \leq RH \qquad (38)$$

is SOC representable where $PM(X,Y)$ is one of the performance metrics: $L^q, W^q, L, W, TW^q, TW$, presented by (2)-(7), and $RH$ can be a constant, variable, or any affine function of the problem's decision variables ($X$ and $Y$ are vectors including all variables $X_p$, $p \in P$ and $Y_{p,r}, p \in P, r \in R$, respectively).

We only deliberate the procedure for two cases $W$ and $TW$, as the other cases can be managed by slight modifications (recall that $TW^q = L^q$, $TW = L = L^q + \lambda\mu^{-1}$, and $W = W^q + \mu^{-1}$). Our unified results are now given in the following two theorems.

**Theorem 1.** Let $PM(X,Y)$ be the expected number of people in the system $L$, (3), or, equivalently, the expected total waiting times in the system $TW$, (7). Then, considering Constraints (34)-(37), Constraint (38) can be represented by (39)-(41),

$$\sum_{p \in P} \frac{(1 + \mu_p^2 \sigma_p^2) r_p}{2\mu_p} + \frac{\sum_{r \in R} \lambda_r Y_{p,r}}{\mu_p} \leq RH \qquad (39)$$

$$\left(\sum_{r \in R} \lambda_r Y_{p,r}\right)^2 \leq r_p \left(\mu_p - \sum_{r \in R} \lambda_r Y_{p,r}\right) \qquad p \in P \qquad (40)$$

$$r_p \geq 0 \qquad p \in P \qquad (41)$$

or equivalently by (42)-(44),

$$\sum_{p \in P} \frac{(1 + \mu_p^2 \sigma_p^2) s_p}{2\mu_p} + \frac{(1 - \mu_p^2 \sigma_p^2) \sum_{r \in R} \lambda_r Y_{p,r}}{2\mu_p} \leq RH \qquad (42)$$

$$\left(\sum_{r\in R}\lambda_r Y_{p,r}\right)^2 \leq \left(s_p\mu_p - \sum_{r\in R}\lambda_r Y_{p,r}\right)\left(\mu_p - \sum_{r\in R}\lambda_r Y_{p,r}\right) \quad p \in P \quad (43)$$

$$s_p \geq 0 \quad\quad p \in P, q \in Q. \quad (44)$$

Moreover, in the case that $Y_{p,r}, p \in P, r \in R$ are binary variables, one can also use (45)-(47) for reformulating (38) as

$$\sum_{p\in P}\frac{(1+\mu_p^2\sigma_p^2)s_p}{2\mu_p} + \frac{(1-\mu_p^2\sigma_p^2)\sum_{r\in R}\lambda_r Y_{p,r}}{2\mu_p} \leq RH \quad (45)$$

$$\sum_{r\in R}\lambda_r Y_{p,r}^2 \leq s_p\left(\mu_p - \sum_{q\in Q}\lambda_q Z_{p,q}\right) \quad p \in P \quad (46)$$

$$s_p \geq 0 \quad\quad p \in P. \quad (47)$$

**Proof.** The proof straightforwardly follows from our arguments given in Section 3, which is not given here for the sake of brevity. ∎

**Remark 1.** We have the following simplifications:
- If the queue is M/M/1, then $\sigma_p^2 = 1/\mu_p^2$ and Constraints (39) and (42) (or (45)) can be simplified as

$$\sum_{p\in P}\frac{r_p + \sum_{r\in R}\lambda_r Y_{p,r}}{\mu_p} \leq RH$$

$$\sum_{p\in P}\frac{s_p}{\mu_p} \leq RH.$$

- If the queue is M/D/1, then $\sigma_p^2 = 0$ and Constraints (39) and (42) (or (45)) are given by

$$\sum_{p\in P}\frac{r_p + 2\sum_{r\in R}\lambda_r Y_{p,r}}{2\mu_p} \leq RH$$

$$\sum_{p\in P}\frac{s_p + \sum_{r\in R}\lambda_r Y_{p,r}}{2\mu_p} \leq RH.$$

**Theorem 2.** Let $PM(X,Y)$ be the expected waiting time in the system $W$, (5). Then, considering (34)-(37), Constraint (38) can be represented by (48)-(50) given below

$$\sum_{p \in P} \frac{(1 + \mu_p^2 \sigma_p^2) s_p}{2\mu_p} + \frac{X_p}{\mu_p} \leq RH \qquad (48)$$

$$\left(\sum_{r \in R} \lambda_r Y_{p,r}\right)^2 \leq \left(s_p \mu_p - \sum_{r \in R} \lambda_r Y_{p,r}\right)\left(\mu_p - \sum_{r \in R} \lambda_r Y_{p,r}\right) \qquad p \in P \qquad (49)$$

$$s_p \geq 0 \qquad p \in P. \qquad (50)$$

Moreover, in the case that $Y_{p,r}, p \in P, r \in R$ are binary variables one can alternatively use (51)-(53):

$$\sum_{p \in P} \frac{(1 + \mu_p^2 \sigma_p^2) s_p}{2\mu_p} + \frac{X_p}{\mu_p} \leq RH \qquad (51)$$

$$\sum_{r \in R} \lambda_r Y_{p,r}^2 \leq s_p \left(\mu_p - \sum_{r \in R} \lambda_q Y_{p,r}\right) \qquad p \in P \qquad (52)$$

$$s_p \geq 0 \qquad p \in P. \qquad (53)$$

**Proof.** Let us define $s_p \geq \frac{\sum_{r \in R} \lambda_r Y_{p,r}}{\mu_p - \sum_{r \in R} \lambda_r Y_{p,r}}$, $p \in P$, then, considering (5), Constraint (38) can be represented as follows:

$$\sum_{p \in P} \frac{(1 + \mu_p^2 \sigma_p^2) s_p}{2\mu_p} + \frac{X_p}{\mu_p} \leq RH.$$

In the case that $Y_{p,r}, p \in P, r \in R$, are binary variables, the constraint $\sum_{r \in R} \lambda_r Y_{p,r} \leq s_p(\mu_p - \sum_{r \in R} \lambda_q Y_{p,r})$ is equivalent to the hyperbolic Constraint (52). However, for the general case where variables $Y_{p,r}, p \in P, r \in R$ are real-valued, one can show that the constraint can be rewritten as Constraint (49) by multiplying both sides by $\mu_p$ and adding $\left(\sum_{r \in R} \lambda_r Y_{p,r}\right)^2$. This completes the proof. ∎

**Remark 2.** Consider the following special cases:

- If the queue is M/M/1, then (48) (or (51)) is given by

$$\sum_{p \in P} \frac{s_p + X_p}{\mu_p} \leq RH.$$

- If the queue is M/D/1, then (48) (or (51)) can be written as

$$\sum_{p \in P} \frac{s_p + 2X_p}{2\mu_p} \leq RH.$$

The above theorems show that for different queue metrics $PM(X,Y)$, Constraint (38) can be reformulated as a set of hyperbolic constrains, which are SOC representable. One can express these constraints in different primary or secondary forms. One should carefully note that an MISOCP solver may be very sensitive to the secondary-form reformulations.

For our specific decision problem when the MISOCP solver is CPLEX, our numerical experiments show that type-I secondary-form constraints mostly perform better. Moreover, the performance of different type-I secondary-form constraints may significantly differ (see 3-MISOCP and 4-MISOCP).

In Remarks 1 and 2, it is shown that the above formulations can be simplified for specific queues such as M/D/1 and M/M/1, but these simplifications do not necessarily improve the computational performance, as it was observed for our location problem.

## 7  Outlook

An important step towards publicizing MISOCP solvers is to reveal how different problems can be molded by MISOCPs. This paper shows that nonlinear M/G/1 queueing formulas incorporated into a decision problem can be cast as MISOCPs when the service rate is represented by a discrete variable and the demand rate depends affinely on some variables. The advantage of using this method is completely demonstrated for a general stochastic location problem in a congested service network.

In future studies, the MISOCP reformulation method can be used in other application areas where the congestion arising in M/D/1, M/M/1 or M/G/1 queues should be controlled. Another open area is investigating the case that the service-capacity is considered an arbitrary decision variable. In this case, the queueing formulas become very complex as the variance term typically depends on the service rate nonlinearly. Extending similar results for G/G/1 queues may be a challenging future research because there are only approximate closed-form formulas available for their performance metrics, and because we no longer have the superposition property of Poisson processes used here to simply analyze the aggregated arrival demand process in M/G/1 queues.

Table 2. Comparison of MISOCPs on test problems used by Vidyarthi and Jayaswal (2014) with 25 facilities, 5 service-capacity levels, and 400 customers; run-time limit is set to 3 hours.

| Test problem | Coefficient of variation* | Average waiting cost per time unit | 1-MISOCP | | | 2-MISOCP | | | 3-MISOCP | | | 4-MISOCP | | |
|---|---|---|---|---|---|---|---|---|---|---|---|---|---|---|
| | | | NO. of nodes | Run time (second), Best gap (%) | Integer-relaxation bound | NO. of nodes | Run time (second), Best gap (%) | Integer-relaxation bound | NO. of nodes | Run time (second), Best gap (%) | Integer-relaxation bound | NO. of nodes | Run time (second), Best gap (%) | Integer-relaxation bound |
| I01 | 0.5 | 1 | <19k | 571, 0.00 | 132057 | <44k | 1665, 0.00 | 117474 | <99k | 3457, 0.00 | 132041 | <16k | 513, 0.00 | 132057 |
| I02 | 0.5 | 50 | <10k | 232, 0.00 | 137982 | <206k | 10800, 0.20 | 121442 | <215k | 10800, 3.78 | 137818 | <7k | 235, 0.00 | 137982 |
| I03 | 0.5 | 500 | <4k | 191, 0.00 | 154740 | <2k | 10800, 1.70 | 134723 | <184k | 10800, 10.08 | 153670 | <2k | 133, 0.00 | 154740 |
| I04 | 0.5 | 5000 | <1k | 52, 0.00 | 229005 | <1k | 10800, 7.26 | 191406 | <196k | 10800, 19.42 | 214691 | <1k | 47, 0.00 | 229005 |
| I05 | 1.5 | 1 | <24k | 530, 0.00 | 132638 | <50k | 1545, 0.00 | 117772 | <237k | 10800, 0.36 | 132612 | <38k | 920, 0.00 | 132638 |
| I06 | 1.5 | 50 | <5k | 272, 0.00 | 141465 | <1k | 10800, 0.99 | 123094 | <239k | 10800, 4.82 | 141102 | <3k | 180, 0.00 | 141465 |
| I07 | 1.5 | 500 | <5k | 246, 0.00 | 165181 | <2k | 10800, 2.10 | 137686 | <203k | 10800, 18.75 | 162379 | <7k | 296, 0.00 | 165181 |
| I08 | 1.5 | 5000 | <5k | 82, 0.00 | 267109 | <1k | 10800, 27.80 | 200954 | <222k | 10800, 31.07 | 211431 | <5k | 70, 0.00 | 267109 |
| I09 | 2.5 | 1 | <15k | 373, 0.00 | 133368 | <94k | 2665, 0.00 | 118176 | <315k | 10800, 0.61 | 133329 | <30k | 715, 0.00 | 133368 |
| I10 | 2.5 | 50 | <4k | 206, 0.00 | 146152 | <2k | 10800, 0.99 | 125105 | <202k | 10800, 7.76 | 145405 | <7k | 273, 0.00 | 146152 |
| I11 | 2.5 | 500 | <5k | 205, 0.00 | 179103 | <1k | 10800, 4.55 | 141338 | <204k | 10800, 22.55 | 172030 | <5k | 170, 0.00 | 179103 |
| I12 | 2.5 | 5000 | <18k | 127, 0.00 | 321830 | <1k | 10800, 80.78 | 215537 | <182k | 10800, 44.35 | 116822 | <32k | 164, 0.00 | 321830 |

* Coefficient of variation: $\mu \times \sigma$
The best run time is highlighted.

Table 3. Comparison of 1-MISOCP, 4-MISOCP, and the existing cutting plane method on test problems used by Vidyarthi and Jayaswal (2014), with 25 facilities, 5 service-capacity levels, and 400 customers; run-time limit is set to 3 hours.

| Test problem | coefficient of variation | Average waiting cost per time unit | Total cost | Establishing cost % | Waiting cost % | Traveling cost % | 1-MISOCP Run time (second), Best gap (%) | 4-MISOCP Run time (second), Best gap (%) | Cutting-plane method Run time (second), Best gap (%) (NO. of iterations) |
|---|---|---|---|---|---|---|---|---|---|
| I01 | 0.5 | 1 | 132360 | 38.57 | 0.36 | 61.07 | 571, 0.00 | 513, 0.00 | 324, 0.00 (2) |
| I02 | 0.5 | 50 | 138262 | 37.10 | 2.80 | 60.10 | 232, 0.00 | 235, 0.00 | 196, 0.00 (3) |
| I03 | 0.5 | 500 | 154964 | 37.99 | 8.38 | 53.63 | 191, 0.00 | 133, 0.00 | 156, 0.00 (3) |
| I04 | 0.5 | 5000 | 229095 | 34.98 | 30.13 | 34.89 | 52, 0.00 | 47, 0.00 | 89, 0.00 (3) |
| I05 | 1.5 | 1 | 132955 | 38.75 | 0.64 | 60.61 | 530, 0.00 | 920, 0.00 | 146, 0.00 (2) |
| I06 | 1.5 | 50 | 141662 | 37.73 | 3.84 | 58.43 | 272, 0.00 | 180, 0.00 | 172, 0.00 (2) |
| I07 | 1.5 | 500 | 165672 | 39.19 | 10.76 | 50.05 | 246, 0.00 | 296, 0.00 | 174, 0.00 (3) |
| I08 | 1.5 | 5000 | 267439 | 35.62 | 34.77 | 29.61 | 82, 0.00 | 70, 0.00 | 424, 0.00 (5) |
| I09 | 2.5 | 1 | 133700 | 39.05 | 0.83 | 60.12 | 373, 0.00 | 715, 0.00 | 385, 0.00 (3) |
| I10 | 2.5 | 50 | 146409 | 38.08 | 5.26 | 56.66 | 206, 0.00 | 273, 0.00 | 238, 0.00 (3) |
| I11 | 2.5 | 500 | 179565 | 41.99 | 14.16 | 43.86 | 205, 0.00 | 170, 0.00 | 239, 0.00 (3) |
| I12 | 2.5 | 5000 | 322098 | 35.48 | 40.40 | 24.11 | 127, 0.00 | 164, 0.00 | 10800, 0.25 (1) |

The best run time is highlighted.

Table 4. Computational results for test problems constructed based on Holmberg et al. (1999), with 30 facilities, 10 service-capacity levels, and 200 customers; run-time limit is set to 3 hours.

| Test problem | Total cost | Establishing cost % | Waiting cost % | Traveling cost % | 1-MISOCP Run time (second), Best gap (%) | 4-MISOCP Run time (second), Best gap (%) | Cutting-plane method Run time (second), Best gap (%) (NO. of iterations) |
|---|---|---|---|---|---|---|---|
| P56 | 2338409 | 47.53 | 16.20 | 27.27 | 6573, 0.00 | 5570, 0.00 | 10800, 3.02 (1) |
| P57 | 2766230 | 48.50 | 12.39 | 39.11 | 10800, 3.09 | 10800, 2.92 | 10800, 6.22 (1) |
| P58 | 3307104 | 45.34 | 6.89 | 47.78 | 10800, 3.26 | 10800, 2.75 | 10800, 3.58 (1) |
| P59 | 2690092 | 42.08 | 9.50 | 48.43 | 639, 0.00 | 646, 0.00 | 10800, 0.04 (2) |
| P60 | 2037424 | 40.44 | 16.71 | 42.85 | 549, 0.00 | 526, 0.00 | 10800, 0.83 (1) |
| P61 | 2423741 | 47.43 | 16.57 | 35.99 | 2541, 0.00 | 9147, 0.00 | 10800, 2.88 (1) |
| P62 | 3107270 | 48.92 | 9.99 | 41.10 | 10800, 2.96 | 10800, 3.19 | 10800, 7.04 (1) |
| P63 | 2408120 | 39.42 | 15.64 | 44.94 | 438, 0.00 | 486, 0.00 | 10800, 0.02 (3) |
| P64 | 1834676 | 37.00 | 18.03 | 44.97 | 240, 0.00 | 291, 0.00 | 10019, 0.00 (4) |
| P65 | 2173829 | 44.16 | 15.73 | 40.11 | 533, 0.00 | 610, 0.00 | 10800, 3.02 (2) |
| P66 | 2856688 | 47.83 | 12.44 | 39.73 | 10800, 2.46 | 10800, 2.85 | 10800, 4.50 (1) |
| P67 | 2429892 | 41.28 | 14.99 | 43.73 | 509, 0.00 | 440, 0.00 | 10800, 0.01 (2) |
| P68 | 1935118 | 38.03 | 16.23 | 45.74 | 273, 0.00 | 247, 0.00 | 2962, 0.00 (4) |
| P69 | 2279064 | 41.61 | 16.25 | 42.14 | 459, 0.00 | 440, 0.00 | 5249, 0.00 (4) |
| P70 | 2908670 | 50.23 | 12.18 | 37.5 | 755, 0.00 | 835, 0.00 | 10800, 0.06 (2) |
| P71 | 2604820 | 41.82 | 13.04 | 45.14 | 4776, 0.00 | 3608, 0.00 | 10800, 4.73 (1) |

The best run time (or the best gap) is highlighted.

Table 5. Comparison of 1-MISOCP and 4-MISOCP on test problems that are constructed by combining each pair of two consecutive test problems given in Table 4, with 60 facilities, 10 service-capacity levels, and 400 customers.

| Test problem | Total cost | Establishing cost % | Waiting cost % | Traveling cost % | 1-MISOCP 5-hour time limit | | 4-MISOCP 5-hour time limit | | 1-MISOCP 10-hour time limit | | 4-MISOCP 10-hour time limit | |
|---|---|---|---|---|---|---|---|---|---|---|---|---|
| | | | | | NO. of nodes | Run time (second), Best gap (%) | NO. of nodes | Run time (second), Best gap (%) | NO. of nodes | Run time (second), Best gap (%) | NO. of nodes | Run time (second), Best gap (%) |
| P56 + P57 | 4143267 | 45.30 | 12.90 | 41.80 | -- | 18000, - | -- | 18000, - | <4k | 29650, 0.00 | <2k | 30593, 0.00 |
| P58 + P59 | 4793767 | 44.35 | 10.90 | 44.75 | -- | 18000, - | -- | 18000, - | <1k | 18330, 0.00 | <1k | 24563, 0.00 |
| P60 + P61 | 3653233 | 44.04 | 17.36 | 38.60 | <2k | 9029, 0.00 | <1k | 9893, 0.00 | <2k | 9029, 0.00 | <1k | 9893, 0.00 |
| P62 + P63 | 4375508 | 39.94 | 11.10 | 48.96 | <1k | 9870, 0.00 | <1k | 10103, 0.00 | <1k | 9870, 0.00 | <1k | 10103, 0.00 |
| P64 + P65 | 3283009 | 40.47 | 19.08 | 40.45 | <11k | 12181, 0.00 | <9k | 13349, 0.00 | <11k | 12181, 0.00 | <9k | 13349, 0.00 |
| P66 + P67 | 4382131 | 44.79 | 14.81 | 40.39 | <2k | 13324, 0.00 | -- | 18000, - | <2k | 13324, 0.00 | <2k | 23129, 0.00 |
| P68 + P69 | 3502339 | 42.08 | 17.84 | 40.08 | <5k | 7891, 0.00 | <9k | 10261, 0.00 | <5k | 7891, 0.00 | <9k | 10261, 0.00 |
| P70 + P71* | 4664334 | 47.16 | 11.93 | 40.92 | 0 | 18000, 1.68 | 0 | 18000, 20.25 | <21k | 36000, 1.32 | <7k | 36000, 1.80 |

--: No feasible solution is found
The best run time (or the best gap) is highlighted.